\documentclass[11pt]{article}       
\usepackage{array}
\usepackage{amsmath}
\usepackage{amssymb}
\begin{document}

\newcommand{\comment}[1]{}    
\newcommand{\hs}{\enspace}
\newcommand{\hhs}{\thinspace}
\newcommand{\real}{\ifmmode {\rm R} \else ${\rm R}$ \fi}
\def\nat{\hbox{\vrule height 7pt width .7pt depth 0pt\hskip -.5pt\bf N}}
\newcommand{\qed}{\hfill{\setlength{\fboxsep}{0pt}
                  \framebox[7pt]{\rule{0pt}{7pt}}} \newline}

\newcommand{\eqed}{\qquad{\setlength{\fboxsep}{0pt}
                  \framebox[7pt]{\rule{0pt}{7pt}}}\newline }
\newtheorem{theorem}{Theorem}
\newtheorem{lemma}[theorem]{Lemma}         
\newtheorem{corollary}[theorem]{Corollary}
\newtheorem{definition}[theorem]{Definition}
\newtheorem{claim}[theorem]{Claim}%
\newtheorem{conjecture}[theorem]{Conjecture}
\newtheorem{proposition}[theorem]{Proposition}
\newtheorem{construction}[theorem]{Construction}
\newtheorem{problem}[theorem]{Problem}

\newcommand{\proof }{{\bf Proof: }}          



\newcommand{\seq}{d_1, \ldots, d_{2k+2}}
\def\eps{\varepsilon}
\newcommand{\eee}{{\mathbb E}}
\newcommand{\seqn}{d_1, \ldots, d_n}
\newcommand{\seqg}{d_1\geq \cdots  \geq d_n}
\newcommand{\edge}{\leftrightarrow}
\newcommand{\nedge}{\nleftrightarrow}
\newcommand{\pakg}{potentially $A_k$-graphical}
\newcommand{\ppkg}{potentially $P_k$-graphical}
\newcommand{\leqk}{\leq _k}
\def\n{\notag}
\def\noedge{\not\leftrightarrow}
\def\f{f_d^k(G)}
\def\fn{f_d^k(K_n)}
\def\fab{f_d^k(K_{a,b})}
\def\to{\rightarrow}
\def\lf{\lfloor}
\def\rf{\rfloor}
\def\lc{\lceil}
\def\rc{\rceil}
\def\Dt{\Delta}
\def\ds{\displaystyle}
\def\tr{\tilde r}
\def\FF{\mathcal{F}}
\def\F{\mathcal{F}}
\def\P{\mathcal{P}}
\def\PP{\mathcal{P}}
\def\GG{\mathcal{G}}
\def\AA{\mathcal{A}}
\def\BB{\mathcal{B}}
\def\HH{\mathcal{H}}
\def\CC{\mathcal{C}}
\def\Q{\mathcal{Q}}
\def\QQ{\mathcal{Q}}
\def\e{\varepsilon}

\title{ Books versus Triangles}
\author{ Dhruv Mubayi
\thanks{ Department of Mathematics, Statistics, and Computer
Science, University of Illinois,  Chicago, IL 60607. Research
supported in part by NSF grant DMS-0653946. Email: {\tt mubayi@math.uic.edu}}\\
}

\date{\today}
\maketitle

\begin{abstract}
A book of size $b$ in a graph is an edge that lies in $b$ triangles. Consider a graph $G$ with  $n$ vertices and $\lfloor n^2/4\rfloor +1$ edges. Rademacher proved that $G$ contains at least $\lfloor n/2\rfloor$ triangles, and Erd\H os conjectured and Edwards proved that $G$ contains a book of size at least $n/6$.

 We prove the following ``linear combination" of these two results. Suppose that $\alpha\in (1/2, 1)$ and the maximum size of a book in $G$ is less than $\alpha n/2$.  Then $G$ contains at least
 $$\alpha(1-\alpha) \frac{n^2}{4}-o(n^2)$$
  triangles as $n\rightarrow \infty$. This is asymptotically sharp.  On the other hand, for every $\alpha\in (1/3, 1/2)$, there exists $\beta>0$ such that $G$ contains at least $\beta n^3$ triangles. It remains an open problem to determine the largest possible $\beta$ in terms of $\alpha$. Our proof uses  the Ruzsa-Szemer\'edi theorem.
\end{abstract}

\section{Introduction}

A book in a graph is a collection of triangles sharing a common
edge. The size of a book is the number of triangles.  Let $b(G)$ be
the size of the largest book in graph $G$ and $t(G)$ be the number
of triangles in $G$. Throughout this note, unless otherwise specified, we let  $G$ be a graph with $n$ vertices and $\lfloor n^2/4\rfloor +1$ edges. All asymptotic notation is to be taken as $n$
grows.

Mantel's theorem states that $G$ contains a triangle, i.e. $t(G) \ge 1$. Rademacher (unpublished) proved in the 1950's that in fact $t(G) \ge \lfloor n/2\rfloor$.
Erd\"os conjectured \cite{E1} in 1962 that $b(G) > n/6$ and this was proved soon after by Edwards (unpublished,  see also Khad\'ziivanov and Nikiforov \cite{KN} for an independent proof).
 Both Rademacher's and Edwards' results are sharp. In the former, $t(G) =\lfloor n/2\rfloor$ is achieved by adding an edge to one part in the complete balanced bipartite graph (note that this also yields $b(G)=\lfloor n/2\rfloor$).  In the latter, every known construction achieving $b(G) = \lfloor n/6\rfloor +1$  has $t(G)=\Omega(n^3)$.

 In this note, we study the relationship between $t(G)$ and $b(G)$. Intuitively, one would suspect that as $t(G)$ decreases, so does $b(G)$. However, this naive intuition is false.
As $t(G)$ becomes smaller, this places greater restrictions on $G$  and $b(G)$ becomes larger,
approaching $n/2$. Indeed, when $t(G)$ is minimized, we saw in the construction above that $b(G)=\lfloor n/2 \rfloor$ which is much larger than $n/6$. On the other hand, when $b(G)=\lfloor n/6\rfloor +1$, which is as small as possible, then $t(G)=\Omega(n^3)$, which is much larger than $\lfloor n/2 \rfloor$.

Our first result shows that as $b(G)$ decreases from $\lfloor n/2\rfloor$ to  $(1-\gamma)n/2$, the number of triangles increases from $\lfloor n/2\rfloor$ to $\Omega_{\gamma}(n^2)$.

\begin{theorem} \label{1} Fix $\alpha \in (1/2, 1)$ and $\e>0$.  Then there exists $n_0$ such that the following holds for  $n>n_0$: Every $n$ vertex graph $G$ with at least $\lfloor n^2/4\rfloor+1$ edges and $b(G)<\alpha n/2$ satisfies
 $$t(G) > (\alpha(1-\alpha)-\e) \frac{n^2}{4}.$$
\end{theorem}

Theorem \ref{1} is asymptotically sharp, as there are examples of graphs
with $t(G)=\alpha(1-\alpha)n^2/4-o(n^2)$ and $b(G) <\alpha n/2$. Indeed,
take the balanced complete bipartite graph (for $n$ even) with one vertex removed and make
this vertex adjacent to $\lfloor \alpha n/2\rfloor-1$ vertices in
one part $n/2-\lfloor \alpha n/2\rfloor +2$ vertices in the
other part.

Our second result shows that if $b(G)<(1/2-\gamma)n/2$, then the number of triangles increases from $\Theta(n^2)$ to $\Omega_{\gamma}(n^3)$.

\begin{theorem}\label{2} For every $\alpha \in (1/3, 1/2)$, there exists
$\beta>0$ such that the following holds for all sufficiently large $n$:
Every $n$ vertex graph $G$ with at least $\lfloor n^2/4\rfloor+1$ edges and $b(G)<\alpha n/2$ has at least
$\beta n^3$ triangles.
\end{theorem}

Note that Theorems \ref{1} and \ref{2} cover all ranges of $\alpha$ except for $\alpha=1+o(1), 1/2+o(1)$, and $1/3+o(1)$.  In particular, $\alpha<1/3$ is impossible due to Edwards' theorem. 

It seems likely that Theorem 2 can be strengthened by replacing
$\beta$ by an explicitly defined  number (in terms of $\alpha$) that is optimal,  but this seems very
hard. One plausible conjecture  is that for every   $\alpha \in (1/3, 1/2)$ one can take
$$\beta={\alpha(1-\alpha)^2
\over 16}+o(1).$$
If true, this would be sharp due to the following
natural generalization of the example achieving equality in Edwards'
theorem:  Partition the vertex set into two almost equal parts $X,
Y$, and partition each of $X$ and $Y$ into three parts roughly of
sizes $\alpha n/2, (1-\alpha)n/4, (1-\alpha)n/4$.  Call them
$X_1, X_2, X_3$ and $Y_1, Y_2, Y_3$ respectively. Add all edges
between the three parts within $X$ and within $Y$. Finally add all
edges between $X_i$ and $Y_i$ for each $i$. 

\section{Tools}
We need the following two results in our proof.  The first is a very special case of the Erd\H os-Simonovits stability theorem \cite{S}.  The proof, which we include here for convenience, is inspired by a recent approach of F\"uredi.
We write $e(G)$ for the number of edges in graph $G$.

\begin{lemma} {\bf (Triangle Stability Lemma)}\label{stab} Let $G$ be a triangle-free graph with $n$ vertices and at least $\lfloor n^2/4 \rfloor -k$ edges. Then $G$ has a vertex partition $X \cup Y$ such that $e(G[X])+e(G[Y]) \le k$.
\end{lemma}

\proof Let $v$ be a vertex of maximum degree.
Since $G$ is triangle-free, there are no edges in $Y:=N(v)$. Let $X=V(G)-N(v)$ and consider the partition $X \cup Y$ of $V(G)$. Let us change $G$ as follows: for each vertex $w \in X$, delete all $s$ edges incident to $w$ contained in $X$ and add $s$ edges from $w$ to $Y$ that were not previously in $G$.  Since $d(w)\le d(v)=|Y|$ this is always possible. Let $G'$ be the graph that results. Now suppose that $G$ has $t$ edges within $X$. Then $e(G')=e(G)+t$  as for every deleted edge within $X$, we add two new edges between $X$ and $Y$. Since $G'$ is bipartite, we have
$$\lfloor n^2/4 \rfloor -k +t=e(G)+t=e(G')\le \lfloor n^2/4 \rfloor.$$
Consequently, $t\le k$ as desired.  \qed

Our second tool is the triangle removal lemma, first proved by Ruzsa and Szemer\'edi \cite{RS}. It is an easy consequence of the Regularity Lemma.

\begin{lemma} {\bf (Triangle Removal Lemma  \cite{RS})} \label{rs}
For every $\delta>0$ there exists $\beta>0$ and $n_0$ such that the following holds for all $n>n_0$:  Every $n$ vertex graph with at most $\beta n^3$ triangles can be made triangle-free by deleting a set of at most $\delta n^2$ edges.
\end{lemma}

\section{Proofs}

In this section we prove Theorems \ref{1} and \ref{2}. Crucial to our proof of Theorem 1 is an assumption on minimum degree, so the theorem that we actually prove is the following:

\begin{theorem} \label{1'}
Let $\alpha'\in (1/2, 1)$.  For every $\e'\in (0, (1-\alpha)/3))$, there exists, $\delta >0$ and $n_0'$ such that the following holds for all $n>n_0'$:
Every $n$ vertex graph $G$ with at least $\lfloor n^2/4\rfloor+1$ edges, minimum degree at least $(1-\delta)n/2$  and $b(G)<\alpha' n/2$ satisfies $t(G) > (\alpha'(1-\alpha')-4\e')n^2/4$.\end{theorem}

Before proceeding to the proof of Theorem \ref{1'} let us argue that it implies Theorem \ref{1}.  We will need the Erd\H os-Stone theorem, which states that for every  $\gamma>0$ there exists $\delta'>0$  and $n_1$ such that every graph with $n>n_1$ vertices and at least $n^2/4+\gamma n^2$ edges contains at least $\delta' n^3$ triangles.
\medskip

{\bf Proof of (Theorem \ref{1'} $\rightarrow$ Theorem  \ref{1}).}
Let us take inputs $\alpha \in (1/2, 1)$ and $\e>0$ from Theorem \ref{1}. Now choose
$$\e' < \min\left \{\frac{\e}{10}\,, \frac{1-\alpha}{2}\right\}.$$
Let $\alpha'=\alpha+\e'$.  The choice of $\e'$ ensures that $\alpha' \in (1/2, 1)$ and $\e'<(1-\alpha')/3$. Let $\delta, n_0'$ be the outputs of Theorem \ref{1'} with inputs $\alpha', \e'$.
We may assume that $\delta <\e'/2$.
Let $\delta'$ and $n_1$ be the outputs of the Erd\H os-Stone theorem with input $\gamma=\delta^2/3$.  Let $n_2$ be sufficiently large so that $\delta'n^3>4n^2$ for all $n>n_2$. 
Finally, let $n_0>2\max\{n'_0, n_1, n_2\}$.

 Now suppose that $n>n_0$ and $G$ is an $n$ vertex graph
with $e(G) \ge \lfloor n^2/4\rfloor+1$ and $b(G)<\alpha n/2$. Our goal is to show that $t(G) > (\alpha(1-\alpha)-\e)n^2/4$.
 
 Note that our constants satisfy the hierarchy
$1/n_0 \ll \delta \ll \e' \ll \e, 1-\alpha$.

 If $G$ has minimum degree  $d < (1-\delta)n/2$,  then
remove a vertex of degree less than $d$ to form the graph $G_1$ with
$n-1$ vertices.  Continue removing a vertex of degree less than
$d_i=(1-\delta)(n-i)/2$ in $G_i$ to form the graph $G_{i+1}$ if such
a vertex exists. Then
$$e(G_k) \ge \left\lfloor {n^2\over 4}\right\rfloor +1 -{(1-\delta)\over 2}\sum_{i=0}^{k-1}(n-i)\ge  {n^2\over 4}-{(1-\delta)\over 2}\sum_{i=0}^{k-1}(n-i). \qquad (1)$$ Suppose that
this procedure continues until   $k =\lceil \delta n\rceil <n/2$. Then by (1) and $n>n_0$ we have
$$e(G_k)\ge {(n-k)^2\over 4} +\left({\delta k n\over 2} +{\delta k\over 4}\right) -\left({k \over
4}+{\delta k^2\over 4}\right)>{(n-k)^2\over 4} +{\delta k n\over
2}-{\delta k^2\over 3}>{(n-k)^2\over 4}+{\delta^2(n-k)^2\over 3}.$$
By the Erd\H os-Stone theorem, $G_k$ (and therefore $G$) has at
least $\delta' (n-k)^3$ triangles and by the choice of $n_0$, this is greater than $4(n-k)^2>n^2$ and we are done. Consequently, $k<\delta n$ and we may assume
that this procedure stops at graph $G_l$ with $n-l>(1-\delta)n$
vertices and at least $(n-l)^2/4+1 $ edges (since the expression in
(1) with $k=l$ is always at least this large).  Since $\delta<1/2$, we have $n-l>n_0/2=n_0'$, and the minimum degree of $G_l$ is at least $(1-\delta)(n-l)/2$, we may try and apply Theorem \ref{1'} to $G_l$. 
The inputs of Theorem \ref{1'} are $\alpha'$ and $\e'$. Because $\delta \ll \e'$ we have 
$$b(G_l)\le b(G)< \frac{\alpha n}{2}<\frac{(\alpha+\e')(1-\delta)n}{2}<\frac{\alpha'(n-l)}{2}.$$
Since $\e', \delta \ll \e$, 
$$(\alpha'(1-\alpha')-4\e')(1-2\delta+\delta^2)
>(\alpha(1-\alpha)-5\e')(1-3\delta)>\alpha(1-\alpha)-\e.$$
Therefore  Theorem \ref{1'} implies that
$$t(G) \ge t(G_l)\ge (\alpha'(1-\alpha')-4\e')\frac{(n-l)^2}{4}
>(\alpha(1-\alpha)-\e)\frac{n^2}{4}.$$
This completes the proof. \qed
\medskip

{\bf Proof of Theorem \ref{1'}.} For notational simplicity, let us replace $\alpha', \e', n_0'$ in Theorem \ref{1'} by $\alpha, \e, n_0$.
So we suppose that $\alpha \in (1/2, 1)$ and $\e\in (0, (1-\alpha)/3)$ are given. Let $\delta=\e^2/50$ and $n_0$ be sufficiently large for all inequalities needed in the proof and for an application of Lemma \ref{rs}. Let $n>n_0$.

Suppose for contradiction that $G$ is a graph with $n$ vertices, at least $\lfloor n^2/4\rfloor+1$ edges, minimum degree at least $(1-\delta)n$, $b(G)<\alpha n/2$ and $t(G) \le (\alpha(1-\alpha)-4\e)n^2/4$.
Since $t(G)<n^2$ and $n>n_0$, by Lemma \ref{rs} we may remove less than $\delta n^2$ edges
from $G$ to make it triangle-free. The resulting graph $G'$ has more
  than $n^2/4-\delta n^2$ edges, so by Lemma \ref{stab}, $G'$ has a vertex partition $A,B$ where $e(G[A])+e(G[B])< \delta n^2$. Now consider a vertex
partition of $G$ into $X,Y$ that maximizes the number of
$X,Y$-edges.  Since one possibility is $A,B$, we are guaranteed that
the number of $X,Y$-edges is at least $n^2/4 -\delta n^2$.
Moreover,  both $X$ and $Y$ have size $(1 \pm 3\sqrt{\delta})n/2$,
otherwise we obtain the contradiction
$${n^2\over 4}<e(G) \le |X||Y|+\delta
n^2<(1-3\sqrt\delta)(1+3\sqrt \delta){n^2\over 4} +\delta
n^2=(1-9\delta){n^2\over 4}+\delta n^2<{n^2\over 4}.$$

Now let $B$ be the set of edges of $G$ entirely contained within $X$
or entirely contained within $Y$, i.e., $B=E(G[X]) \cup E(G[Y])$.
Let $M$ be the set of pairs in $X \times Y$ that are not edges of
$G$. Then $E(G)-B \cup M$ is bipartite, so it has at most $n^2/4$
edges. As $e(G)\ge \lfloor n^2/4\rfloor +1$, we conclude that
$$|M|< |B|<\delta n^2.$$
In particular, $B\ne\emptyset$.
Next, let $M' \subset M$ be the set of those pairs $\{x,y\} \in M$,
such that $x$ and $y$ are each incident with at least $\e n$
edges of $B$ (of course $x$ and $y$ are on opposite sides of the
partition).

{\bf Claim.}  $|M'| <c:=\lceil \alpha^2(1-\alpha)^2/\e^4\rceil $.

{\bf Proof of Claim.} Otherwise, by the K\"onig-Hall theorem, there
is either a matching or a star of size at least
$s=\lfloor \alpha(1-\alpha)/\e^2\rfloor $ in $M'$. In the case of a matching,
each pair $f=uv$ of the matching is incident with $\lceil \e n\rceil$ edges
of $B$ in  both $X$ and $Y$. Consider any set of $\lceil \e n\rceil$ edges of $B$
incident to $u \in X$. By the choice of $X,Y$, each vertex has at
least as many neighbors on the opposite side of the partition as its
own side, hence $u$ has at least $\e n$ neighbors in $Y$. Each
edge between these two sets of neighbors of $u$ forms a triangle,
and the number of such edges is at least
$$\e^2n^2 -|M|>\e^2n^2-\delta n^2>{\e^2\over 2} n^2.$$ Every two
such pairs $uv, u'v'$ in the matching of $M'$ count at most $4n$
common triangles, so by Inclusion/Exclusion, we obtain the
contradiction
$$t(G) >2s{\e^2\over 2} n^2 - {2s \choose 2} 4n>
s\e^2 n^2-8s^2n\ge\alpha(1-\alpha)n^2-\e^2 n^2-{8\over
\e^4}n>(\alpha(1-\alpha)-\e){n^2\over 4}.$$

In the case of a star, the same argument works, since we count
triangles starting from a vertex in the part corresponding to the
leaf set of this star of $M'$.  Indeed, suppose we have pairs $xy_1,
\ldots, xy_s \in M'$.  Then for each $i$, consider $\lceil \e n\rceil$
edges of $B$ incident with $y_i$. Since $y_i$ can be moved to the
other part, there are at least $\e n$ edges of $E(G)-B$
incident with $y_i$. Now proceed to find many triangles as in the
previous case. \qed

Form the bipartite graph $H$ with parts $B$ and $M$, where $e \in B$
is adjacent to $f \in M$ if edge $e$ is incident with the pair $f$.
Since $b(G) \le \alpha n/2$, every $e \in B$ is adjacent
(in $H$) to at least
$$\min\{|X|, |Y|\}-\frac{\alpha n}{2}>(1-3\sqrt \delta -\alpha)\frac{n}{2}$$
vertices $f \in M$. Consequently
$$e(H)\ge |B|(1-3\sqrt \delta-\alpha){n \over 2}.$$
The number of edges of $H$ incident to $M'$ is at most $|M'||B| < c|B|$
so the number of edges in $H$ incident to some pair of $M-M'$ is at least
$$|B|(1-\alpha-3\sqrt{\delta})\frac{n}{2}- c|B| >
|B|(1-\alpha+\e)\frac{n}{2}.$$ Since $|M| <|B|$, we conclude that
there is an $f=uv \in M-M'$ that is incident (in $G$) with at least
$(1-\alpha-\e)n/2>\e n$ distinct $e \in B$.  Since $f  \not\in
M'$, we may assume (wlog) that at least
$$(1-\alpha-\e)\frac{n}{2}-\e n>(1-\alpha-3\e)\frac{n}{2}$$
of these edges $e$ lie in $X$, say they
form a star in $G$ with center $u$ and leaf set $L_X=N(u) \cap X$. 
So we have $|L_X|>(1-\alpha-3\e)n/2$.
Let $L_Y=N(u) \cap Y$.

If $|L_Y|<|L_X|$, then we could move $u$ to $Y$ and increase the number of edges between $X$ and $Y$, thereby contradicting the choice of the partition $X \cup Y$. We therefore have $|L_Y| \ge |L_X|$. 
As $G$ has
minimum degree at least $(1-\delta)n/2$, we have $d(u)=|L_X|+|L_Y|\ge (1-\delta)n/2$.  Consequently,
$$|L_Y|\ge \max\left\{|L_X|, \frac{(1-\delta)n}{2}-|L_X|\right\}.$$ Let $|L_X|=an/2$ and $|L_Y|=bn/2$. The number of edges in $G$ between $L_X$ and
$L_Y$ is at least
$$|L_X||L_Y|-|M| \ge ab\frac{n^2}{4}-\delta n^2=(ab-4\delta){n^2 \over 4},$$
where
$$b \ge a \ge 1-\alpha-3\e \qquad \hbox{ and }\qquad a+b \ge 1-\delta.$$
Now $ab-4\delta$ is minimized by minimizing $a+b$ and then maximizing $b-a$. Since $\e <(1-\alpha)/3$, the minimum occurs at
 $$a=1-\alpha-3\e>0  \quad \hbox{ and } \quad
b=1-\delta-a=\alpha-\delta+3\e$$ where it equals
$$(1-\alpha-3\e)(\alpha-\delta+3\e)
-4\delta>\alpha(1-\alpha)-4\e.$$ Since each  of these edges
gives rise to a unique triangle, we conclude that $t(G)>(\alpha(1-\alpha)-4\e)n^2/4$, a contradiction.  \qed

\bigskip

{\bf Proof of Theorem 2.} We use the notation from Theorem
\ref{1'}'s proof. Let $\alpha \in (1/3, 1/2)$ be given and choose
$$\e=\min\left\{\frac{1}{10}\,, \frac{1-2\alpha}{4}\right\}.$$
Note that $\alpha<1/2$ implies that $\e>0$.
Let $\delta=\e^2/50$  and let $\beta$ be sufficiently small
so that we can apply Lemma \ref{rs} with input $\delta$ and output $\beta$. Our hierarchy of constants is
$$1/n_0 \ll \beta \ll \delta \ll \e \ll \alpha.$$
We do not need  the minimum degree assumption on $G$.
Suppose for contradiction that
$b(G) <\alpha n/2$ and $t(G)<\beta n^3$. By Lemma \ref{rs}
we can make $G$ triangle-free by removing a set of at most $\delta n^2$ edges. Now follow the proof of Theorem \ref{1'} precisely to obtain the partition $X, Y$ with the same properties and also $|M|<|B|<\delta n^2$. We may also assume the Claim from Theorem \ref{1'}'s proof holds.
 Then we find a pair $f=uv \in
M-M'$ incident with at least $(1-\alpha-\e)n/2$ distinct $e \in
B$.  Again form the sets $L_X$ and $L_Y$ whose vertices are neighbors of $u \in X$. By optimality of the partition, we  have
$$|L_Y|\ge |L_X|\ge (1-\alpha-3\e)\frac{n}{2}>\frac{n}{4}.$$
 Consider the subgraph $K$ of
edges of $G$ between $L_X$ and $L_Y$.  Then
$$e(K)\ge|L_X||L_Y|-|M|$$
so there exists a vertex $v \in L_X$
with
$$d_K(v) \ge \frac{e(K)}{|L_X|}\ge |L_Y|-\frac{|M|}{|L_X|}
\ge (1-\alpha-3\e)\frac{n}{2}-\frac{\delta n^2}{n/4}>
(1-\alpha-4\e)\frac{n}{2}.$$
Since $\alpha <1/2$ and $\e<(1-2\alpha)/4$, this is at least $\alpha n/2$.
Therefore the edge $uv$ lies in at least $d_K(v)\ge \alpha n/2$ triangles, contradicting the hypothesis $b(G)<\alpha n/2$.
\qed
\section{Concluding Remarks}

$\bullet$ We observed that as $\alpha$ decreases from 1 to 1/3, the number of triangles increases from $\lfloor n/2 \rfloor$ to $\Omega(n^3)$.
Theorem 1 shows that for $\alpha<1$ we always have $t(G)=\Omega(n^2)$, Similarly, as $\alpha$ changes from $1/2+o(1)$ to $1/2-o(1)$, the number of triangles changes from quadratic to cubic in $n$.  There appear to be  two phase transitions here, $\alpha=1$ and $\alpha=1/2$.  It would be very interesting to understand the scaling window in these two ranges, namely, the rate at which $t(G)$ changes from linear to quadratic and from quadratic to cubic.  Perhaps this is connected with the behavior of various parameters in the regularity Lemma, since this seems crucial to our argument.
\medskip

$\bullet$  One could also ask the same questions for graphs with $\lfloor n^2/4 \rfloor +q$ edges for $q>1$. Results of Erd\" os \cite{E1} and Lov\'asz-Simonovits \cite{LS} determine the minimum number of triangles  and Bollob\'as and Nikiforov \cite{BN} determine the minimum value of $b(G)$. Theorems \ref{1} and \ref{2} apply to this case, since the hypothesis is simply $e(G) \ge  \lfloor n^2/4 \rfloor +1$.  Moreover, when $q=o(n)$, the results are asymptotically sharp, as evidenced by easy modifications of the constructions shown earlier.
The situation when $q=\Omega(n)$ appears to be more complicated and our methods do not seem to apply.

\medskip

$\bullet$
One could consider cliques of larger size and the appropriately defined books (collection of cliques that share an edge).  Our proofs appear to be robust enough to address this situation in a similar fashion, in particular, the tools we need
(removal lemmas, stability results, results for generalized books) are available. Nevertheless, the technical details would probably be quite complicated, and since the situation for triangles is not yet well understood, we have chosen not to address this.


\begin{thebibliography}{9}

\bibitem{BN} B. Bollob\'as, V. Nikiforov,  Books in graphs,  European J. Combin. 26 (2005), 259-270


\bibitem{E1} P. Erd\H os, On a theorem of Rademacher-Tur\'an, Illinois Journal of Math, 6, (1962), 122--127

\bibitem{KN} N. Khad\'ziivanov, V. Nikiforov, Solution of a problem of P. Erd\H os about the maximum number of triangles with a common edge in a graph (Russian), C. R. Acad Bulgare Sci 32 (1979) 1315--1318


\bibitem{LS} L.   Lov\'asz, M. Simonovits,  On the number of complete subgraphs of a graph. II,  Studies in pure mathematics,  459--495, Birkhäuser, Basel, 1983.

\bibitem{RS}
I. Ruzsa, E. Szemer\'edi, Triple systems with no six points carrying three triangles, Combinatorics (Proc. Fifth Hungarian Colloq., Keszthely, 1076), Vol. II, Colloq. Math. Soc. J\'anos Bolyai, vol. 18, North-Holland, Amsterdam, 1978, pp. 939--945.



\bibitem{S} M. Simonovits,  A method for solving extremal problems in graph theory, stability problems, 1968  Theory of Graphs (Proc. Colloq., Tihany, 1966)  pp. 279--319 Academic Press, New York
\end{thebibliography}
\end{document}